\documentclass[10pt]{article}

\usepackage{graphicx}
\usepackage{subcaption}

\usepackage[a4paper, left=35mm,right=35mm,top=34mm,bottom=34mm]{geometry}
\usepackage[utf8]{inputenc}
\usepackage[T1]{fontenc}
\usepackage[english]{babel}

\usepackage{enumerate}
\usepackage{graphicx}
\usepackage{hyperref}
\hypersetup{
    colorlinks=true,
    linkcolor=blue,
    filecolor=magenta,      
    urlcolor=cyan,
}

\usepackage{mathtools,amsthm,amssymb,amsfonts}
\usepackage{algorithm}
\usepackage{algpseudocode}
\makeatletter
\def\thm@space@setup{
  \thm@preskip=10pt \thm@postskip=10pt
}
\makeatother
\theoremstyle{plain}

\theoremstyle{plain}

\theoremstyle{definition}

\theoremstyle{definition}

\theoremstyle{remark}

\theoremstyle{remark}

\usepackage{caption} 
\usepackage{booktabs}
\usepackage{listings}
\usepackage{color}
\definecolor{dkgreen}{rgb}{0,0.6,0}
\definecolor{gray}{rgb}{0.5,0.5,0.5}
\definecolor{mauve}{rgb}{0.58,0,0.82}
\usepackage{enumitem}

\newcommand{\abs}[1]{\left\lvert#1\right\lvert}
\newcommand{\norm}[1]{\left\lVert#1\right\rVert}

\DeclareMathOperator*{\argmin}{arg\,min}

\newcommand{\email}[1]{\protect\href{mailto:#1}{#1}}

\usepackage{xcolor}[2.11]
\colorlet{inlinkcolor}{green!50!black}
\colorlet{exlinkcolor}{red!50!black}
\if@hidelinks
\hypersetup{ hidelinks = true }
\else
\hypersetup{
  colorlinks = true,
  allcolors = inlinkcolor,
  urlcolor = exlinkcolor,
}
\fi

\newenvironment{@abssec}[1]{
        \vspace{.05in}\parindent .0in
        {\upshape\bfseries #1. }\ignorespaces
    }
    {\par\vspace{.1in}}
\renewenvironment{abstract}{\begin{@abssec}{\abstractname}}{\end{@abssec}}
\newenvironment{keywords}{\begin{@abssec}{Keywords}}{\end{@abssec}}

\usepackage{fancyhdr}

\author{
  {\normalsize Qinmeng Zou}\thanks{CentraleSup\'elec, Universit\'e Paris-Saclay, 3 rue Joliot Curie, 91190 Gif-sur-Yvette, France
    (\email{zouqinmeng@gmail.com}, \email{frederic.magoules@hotmail.com}).}
  \and
  {\normalsize Fr\'ed\'eric Magoul\`es\footnotemark[1]}
}
\title{On Extensions of Limited Memory Steepest Descent Method}
\date{}

\begin{document}
\maketitle
\thispagestyle{fancy}

\begin{abstract}
We present some extensions to the limited memory steepest descent method based on spectral properties and cyclic iterations.
Our aim is to show that it is possible to combine sweep and delayed strategies for improving the performance of gradient methods.
Numerical results are reported which indicate that our new methods are better than the original version.
Some remarks on the stability and parallel implementation are shown in the end.
\end{abstract}


\section{Introduction}

Unconstrained nonlinear optimization is an exciting field and has been the subject of research for centuries.
One of the most widely used method is the gradient-based iterative method
\begin{equation}\label{eq:x}
x_{n+1} = x_n - \alpha_n g_n,
\end{equation}
where $g_n=\nabla f(x_n)$ is the gradient vector and $\alpha_n>0$ is viewed as a steplength.
Consider a strongly convex quadratic function $f:\mathbb{R}^N\rightarrow\mathbb{R}$ in the form
\begin{equation}\label{eq:f}
f(x) = \frac{1}{2}x^\intercal Ax - b^\intercal x,
\end{equation}
where $A$ is an $N$-dimensional symmetric positive definite (SPD) matrix and $b$ is an $N$-dimensional vector.
The unique minimizer
\[
\min_{x\in\mathbb{R}^N} f(x)
\]
can be obtained by proper gradient methods.
In the past three decades, gradient methods have regained interest with the work of Barzilai and Borwein~\cite{Barzilai1988}, in which a novel approach, the so-called Barzilai-Borwein method (BB), was proposed in the form
\[
\alpha_n^\text{BB} = \frac{s_{n-1}^\intercal s_{n-1}}{s_{n-1}^\intercal y_{n-1}},
\]
where $s_{n-1}=x_n-x_{n-1}$ and $y_{n-1}=g_n-g_{n-1}$.
This method is nonmonotone in regard to both the sequences~$\{f(x_n)\}$ and~$\{\norm{g_n}\}$ where $\norm{g_n}$ denotes the $2$-norm of vector~$g_n$.
Many generalizations of BB can thus be considered, see~\cite{diSerafino2018} and the references therein.

In 2012, Fletcher~\cite{Fletcher2012} proposed a sweep method called limited memory steepest descent (LMSD).
During a sweep, spectral information is exploited in order to approximate the eigenvalues of the Hessian matrix~$A$.
It follows from~\eqref{eq:x} and~\eqref{eq:f} that
\begin{equation}\label{eq:g}
g_{n+1} = g_n - \alpha_nAg_n.
\end{equation}
Let $\{\lambda_i\}$ denote the set of eigenvalues of~$A$ with $i\in\{1,\,\dots,\,N\}$ and $\{v_i\}$ denote the set of corresponding eigenvectors.
We assume that $\lambda_1 \le \dots \le \lambda_N$ and notice that $g_n$ can de decomposed as
\[
g_n = \sum_{i=1}^N \zeta_{i,n}v_i,
\]
where $\zeta_{i,n}$ is the $i$th spectral component of~$g_n$.
It follows that
\[
\zeta_{i,n+1} = (1-\alpha_n\lambda_i)\zeta_{i,n} = \prod_{j=0}^n(1-\alpha_j\lambda_i)\zeta_{i,0}.
\]
If there exists any $\alpha_j=1/\lambda_i$, then the $i$th spectral component vanishes for all subsequent iterations.
In particular, if $A$ is a diagonal matrix, then $\zeta_{i,n}$ can be denoted by~$g_{i,n}$, expressing the $i$th component of~$g_n$.
Fletcher came up with the idea of using the back gradient vectors to compute the Ritz values, regarded as estimates of the eigenvalues of~$A$.
Consider the most recent $m$ back gradient vectors
\[
G = [g_{n-m},\,g_{n-m+1},\,\dots,\,g_{n-1}].
\]
The QR factorization of~$G$ yields an $N\times m$ matrix~$Q$ with orthonormal columns and an $m\times m$ upper triangular matrix~$R$.
By~\eqref{eq:g}, we know that the columns of~$G$ lie in the Krylov subspace
\[
\text{span}\left\{g_{n-m},\,Ag_{n-m},\,\dots,\,A^{m-1}g_{n-1}\right\}.
\]
Hence, the orthogonalization can be regarded as a Lanczos process (see,~e.g.,~\cite{Magoules2015e}) which leads to a tridiagonal matrix
\begin{equation}\label{eq:t}
T = Q^\intercal AQ.
\end{equation}
The eigenvalues of~$T$, namely, the Ritz values of Hessian matrix~$A$, can be easily solved, for instance by QL algorithm with implicit shift as mentioned in~\cite{Fletcher2012}.

For general unconstrained optimization problems, $A$ is not available.
The following equations provide another way to formulate this method:
\[
T = [R,\,r_n] J R^{-1},\quad J=
\left(\begin{array}{ccc}
\alpha_{n-m}^{-1} & & \\[2pt]
-\alpha_{n-m}^{-1} & \ddots & \\[2pt]
& \ddots & \alpha_{n-1}^{-1} \\[2pt]
& & -\alpha_{n-1}^{-1}
\end{array}\right),
\]
where $R$ and $r$ can be computed through the partially extended Cholesky factorization
\begin{equation}\label{eq:chol}
G^\intercal [G,\,g_n] = R^\intercal [R,\,r].
\end{equation}
It is noteworthy that in the case of minimizing quadratic function~\eqref{eq:f} and with~$m=1$, LMSD is mathematically equivalent to BB.

\section{Extensions of LMSD method}

Extensions to LMSD can be thought of by using techniques in the field of delayed gradient methods.
The first approach is motivated by the spectral properties.
Consider the Yuan steplength~\cite{Yuan2006} in the form
\[
\alpha_n^\text{Y} = \frac{2}{\alpha_n^\text{RA}+\sqrt{\left(\alpha_n^\text{RA}\right)^2 - 4\Gamma_n}},
\]
where
\[
\alpha_n^\text{RA} = \frac{1}{\alpha_{n-1}^\text{SD}}+\frac{1}{\alpha_n^\text{SD}},\quad \Gamma_n = \frac{1}{\alpha_{n-1}^\text{SD}\alpha_n^\text{SD}}-\frac{\norm{g_n}^2}{\left(\alpha_{n-1}^\text{SD}\norm{g_{n-1}}\right)^2}.
\]
Under some assumptions without loss of generality, the following limit holds:
\[
\lim_{n\rightarrow\infty}\alpha_n^\text{Y} = \frac{1}{\lambda_N}.
\]
We refer the reader to~\cite{ZouX2} for more details.
On the other hand, the most basic gradient method is steepest descent (SD) which can be written in the form
\[
\alpha_n^\text{SD} = \argmin_{\alpha}f(x_n-\alpha g_n).
\]
A method that relies on these two stepelengths was proposed in~\cite{DeAsmundis2014} which can be summarized as follows:
\begin{enumerate}
\item Execute several gradient iterations based on SD steplength.
\item Compute $\alpha_n^\text{Y}$ as a constant using data from the current and last iterations.
\item Execute several gradient iterations based on this constant steplength.
\end{enumerate}
For the quadratic case, this method can foster the sequence~$\{1/\alpha_n^\text{Y}\}$ to approximate some largest eigenvalues, and thus reduce the search spaces into smaller and smaller dimensions (see~\cite{DeAsmundis2014}).

Our main idea is to consider the SD iterations being used for both fostering alignment and generating $m$-dimensional subspaces.
Here we report the main steps:
\begin{enumerate}
\item Execute $m$ gradient iterations based on SD steplength.
\item Compute $\alpha_n^\text{Y}$ as a constant using data from the current and last iterations.
\item Execute $d$ gradient iterations based on this constant steplength.
\item Compute $m$ Ritz values of the Hessian matrix.
\item Execute $m$ gradient iterations based on these Ritz values.
\end{enumerate}
We observe that the algorithm outlined above is a combination of alignment method and LMSD.
This might take advantage of both strategies and improvements might be expected.
There exist other constant steplengths that could be exploited for the purpose of alignment, see~\cite{ZouX1} for more details.

Another line of extensions starts with the cyclic steplength
\[
\alpha_n = \alpha_{n-1}.
\]
A theorem provided in~\cite{Friedlander1999} could prove its convergence when $\alpha_n$ follows the framework of gradient method with retards.
It was recently shown in~\cite{Zou2018,ZouX3} that cyclic formulation may accelerate the iteration process both in terms of sequential and parallel computation.
The main steps based on this strategy can be stated as follows:
\begin{enumerate}
\item Compute $m$ Ritz values of the Hessian matrix.
\item Execute $k\cdot m$ gradient iterations based on these Ritz values.
\end{enumerate}
Each Ritz value would be used for $k$ times for updating iterates.
Only $m$ back gradient vectors will be stored for computing the next set of Ritz values.

The previous methods are called LMSD, LMSD with constant steplength (LMSDC), LMSD with retards (LMSDR), respectively.
For LMSD, Fletcher~\cite{Fletcher2012} has further studied two questions: how can the first $m$ back gradient vectors be initialized and in what order should Ritz values be used for performing the gradient iterations within a sweep?
He suggested to initialize one or more Ritz values and proceed in a greedy way, namely, in what follows one could use all back gradient vectors for constructing~$G$ until it is possible to use $m$ back vectors.
For LMSDC, it is clear that SD steplengths must be employed in the first iterations.
For LMSDR, the initialization remains the same as that of LMSD and thus could be carried out by greedy choices as mentioned above.

In~\cite{Fletcher2012}, Ritz values are selected in decreasing order, which gives the best chance that the monotonicity could be retained in early steps.
In addition, if nonmonotonicity is detected in terms of~$f(x_n)$ or~$\norm{g_n}$, then the current sweep is terminated.
Obviously, these approaches can be adapted for LMSDC without additional operations.
For LMSDR, since $m$ Ritz values are employed in $k\cdot m$ iterations, the most reliable approach is to execute all these iterations in decreasing order with the monotonicity detector mentioned above.

\section{Numerical Experiments}

We implement these algorithms by a linear algebra programming library~\cite{Magoules2015}.
QR factorization combining with~\eqref{eq:t} is used for computing the Ritz values.
There are two reasons: the first one is our programming environment has enough working memory that can ensure the storage of several long vectors; the second one is the implementation based on Cholesky factorization suffers from numerical instabilities since $R$ can be numerically ill-conditioned, see~\cite{Fletcher2012} for more details.
Hence, a positive definite quadratic function is used as test case.
We choose the same problem as that in~\cite{Fletcher2012}, which is based on an SPD matrix of size~$20$ with~$\lambda_1=1$.
Other eigenvalues are distributed in geometric progression with ratio~$\sqrt{2}$.
More specifically, the Hessian matrix can be written as follows:
\[
A = \text{diag}(1,\,1.414,\,\dots,\,724.077).
\]
The right hand side~$b$ is selected as a random vector in which the elements range from~$10$ to~$20$.
The stopping criterion is~$\norm{g_n}<10^{-6}\norm{g_0}$ and $x_0$ is fixed with zero.

\begin{figure}[!t]
\centering
\includegraphics[width=4in]{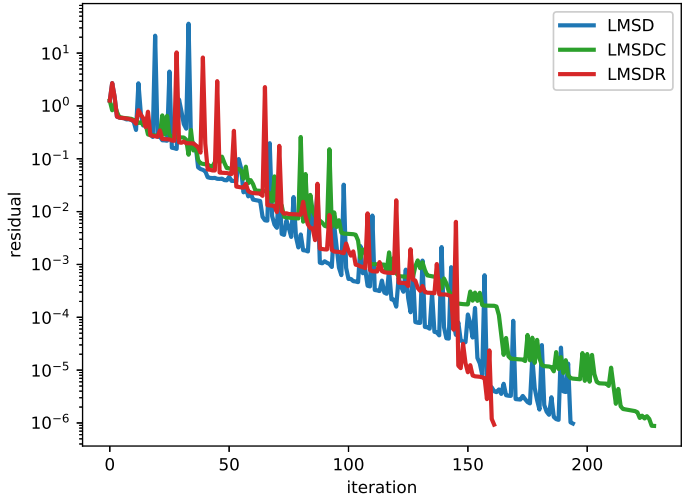}
\caption{Convergence history of gradient norm (residual)~$\norm{g_n}$}
\label{fig:1}
\end{figure}

\begin{figure}[!t]
\centering
\includegraphics[width=4in]{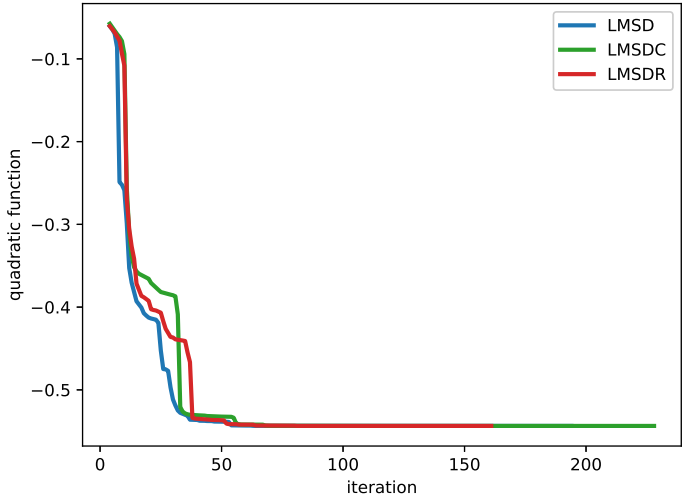}
\caption{Convergence history of quadratic function~$f(x_n)$}
\label{fig:2}
\end{figure}

In Figures~\ref{fig:1} and~\ref{fig:2} we show the history of residual and quadratic function, respectively.
We choose~$m=4$, $d=4$, and $k=2$.
The figures show that the convergence behaviors of LMSD and its two variants are about the same.
Although $\norm{g_n}$ is nonmonotone throughout iterations, $f(x_n)$ is monotonically decreased.
Table~\ref{tab:1} provides computation times, from which we find that LMSDC and LMSDR are significantly better than the original version.
\begin{table}[!t]
\renewcommand{\arraystretch}{1.3}
\caption{Average running time among ten tests and counts of nonmonotone termination within a sweep}
\label{tab:1}
\centering
\begin{tabular}{|c||c|c|c|}
\hline
Method & LMSD & LMSDC & LMSDR \\
\hline
Time (s) & $0.049$ & $0.030$ & $0.031$ \\
$f(x_{n+1}) > f(x_n)$ & $16$ & $1$ & $12$ \\
$\norm{g_{n+1}} > \norm{g_n}$ & $20$ & $10$ & $15$ \\
\hline
\end{tabular}
\end{table}
This observation is consistent with our expectations, which shows better performance of cyclic iterations in terms of computation costs compared with classical iterations.

We give the results of nonmonotone termination within a sweep, shown in Table~\ref{tab:1}.
We can see that LMSDC almost preserves monotonicity even without local termination criterion.
On the other hand, we need to use the termination strategy mentioned above to provide a certain monotonicity for LMSD and LMSDR.

\begin{figure}[!t]
\centering
\includegraphics[width=4in]{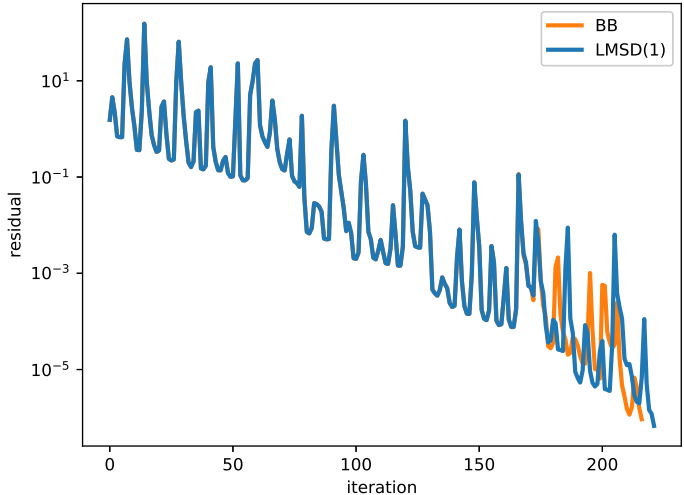}
\caption{Comparison of BB and nonmonotone LMSD with~$m=1$}
\label{fig:3}
\end{figure}

It is known that LMSD is mathematically equivalent to BB for~$m=1$ if nonmonotone sweeps are considered.
However, such behavior may not occur in experiments which is subject to rounding errors.
An example is shown in Figure~\ref{fig:3}.
Note that we use here the same Hessian matrix and right hand side as the preceding test.
In the figure we observe that the convergence results of BB and LMSD are about the same in the beginning.
From about iteration~$170$ we can see some significant difference in colors which indicates the effect of rounding errors.
In Table~\ref{tab:2} we highlight some points in the convergence history and notice that when $n=200$ the absolute difference of residuals even exceeds the second residual value.
This observation is consistent with the two curves in Figure~\ref{fig:3}, which proves that BB and nonmonotone LMSD with~$m=1$ are not numerically equivalent but only mathematically equivalent.
\begin{table}[!t]
\renewcommand{\arraystretch}{1.3}
\caption{Highlight of residuals and their absolute difference for some BB and nonmonotone LMSD(1) iterations}
\label{tab:2}
\centering
\begin{tabular}{|c||c|c|c|c|c|}
\hline
Iteration & $100$ & $150$ & $200$ \\
\hline
BB ($p_1$) & $2.05\times10^{-3}$ & $1.49\times10^{-2}$ & $2.39\times10^{-5}$ \\
LMSD(1) ($p_2$) & $2.05\times10^{-3}$ & $1.49\times10^{-2}$ & $6.45\times10^{-6}$ \\
$\abs{p_1-p_2}$ & $5.57\times10^{-11}$ & $2.12\times10^{-7}$ & $1.75\times10^{-5}$ \\
\hline
\end{tabular}
\end{table}

\section{Practical Considerations}

Experience shows that computing Ritz values according to~\eqref{eq:t} can be a better choice in terms of stability at the expense of more arithmetic operations and extra long vectors.
Fletcher~\cite{Fletcher2012} gave some remedies which consist of discarding the oldest back gradient vectors and recomputing matrix~$T$.
These approaches were also reported in~\cite{diSerafino2018} as parts of the implementation.
Note that cyclic formulation may amplify the rounding errors.
As a result, we could not choose a large~$k$ in practice.
On the other hand, a smooth reduction in residual can be expected for the alignment process, and in fact the local termination detector can be removed without changing monotonicity.

The matrix~$G$ with~$N\gg m$ can be factorized by parallel algorithms.
For example, Golub and Van Loan~\cite{Golub2013} provided some classical strategies for parallel factorizations.
One could employ asynchronous iterations (see~,e.g,~\cite{Magoules2018c,Magoules2018}) in an iterative methods for improving the parallel performance. 
However, it is difficult to handle the effect of dot product operations.
For the communication-avoiding techniques, Demmel et al.~\cite{Demmel2012} discussed parallel QR factorization for tall and skinny matrices which was compared with the classical Householder QR.
The storage of long vectors can be avoided by using Cholesky factorization~\eqref{eq:chol}, which has also been discussed in the previous references.
On the other hand, the communication cost is reduced by cyclic formulation in LMSDR, which may lead to better parallel performance.

We have already implemented the parallel variants by JACK~\cite{Magoules2017b,Magoules2018b}, which is an MPI-based communication library with several choices of convergence detectors~\cite{Magoules2018e} for both classical and asynchronous iterative computing.
Experiments show that LMSDR often gives better convergence results.
There exist other possibilities that one could implement parallel LMSD variants by overlapping communication phases with computations.
Specifically, the bottleneck in these methods is the dot product operations, which can be tackled by synchronization-reducing techniques.
We refer to~\cite{ZouX4} and the references therein for more detailed discussion on this topic.

\section{Conclusions}

We have presented two variants for the LMSD method, called LMSDC and LMSDR, respectively.
The main ingredients that are exploited when formulating the new methods are alignment and cyclic strategies, which have been widely used in delayed gradient methods.
The comparison on the problem used in~\cite{Fletcher2012} reveals that our new methods are superior in terms of computation time and generate similar convergence curves as the original version.

\section*{Acknowledgment}

This work was funded by the project ADOM (M\'ethodes de d\'ecomposition de domaine asynchrones) of the French National Research Agency (ANR).

\bibliography{ref}

\begin{thebibliography}{10}

\bibitem{Barzilai1988}
J.~Barzilai and J.~M. Borwein.
\newblock Two-point step size gradient methods.
\newblock {\em IMA J. Numer. Anal.}, 8(1):141--148, 1988.

\bibitem{DeAsmundis2014}
R.~{De Asmundis}, D.~{di Serafino}, W.~W. Hager, G.~Toraldo, and H.~Zhang.
\newblock An efficient gradient method using the {Y}uan steplength.
\newblock {\em Comput. Optim. Appl.}, 59(3):541--563, 2014.

\bibitem{Demmel2012}
J.~W. Demmel, L.~Grigori, M.~Hoemmen, and J.~Langou.
\newblock Communication-optimal parallel and sequential {QR} and {LU}
  factorizations.
\newblock {\em SIAM J. Sci. Comput.}, 34(1):A206--A239, 2012.

\bibitem{diSerafino2018}
D.~{di Serafino}, V.~Ruggiero, G.~Toraldo, and L.~Zanni.
\newblock On the steplength selection in gradient methods for unconstrained
  optimization.
\newblock {\em Appl. Math. Comput.}, 318:176--195, 2018.

\bibitem{Fletcher2012}
R.~Fletcher.
\newblock A limited memory steepest descent method.
\newblock {\em Math. Program.}, 135(1):413--436, 2012.

\bibitem{Friedlander1999}
A.~Friedlander, J.~M. Mart{\'i}nez, B.~Molina, and M.~Raydan.
\newblock Gradient method with retards and generalizations.
\newblock {\em SIAM J. Numer. Anal.}, 36(1):275--289, 1999.

\bibitem{Golub2013}
G.~H. Golub and C.~F. {Van Loan}.
\newblock {\em Matrix Computations}.
\newblock Johns Hopkins University Press, 4th edition, 2013.

\bibitem{Magoules2015}
F.~Magoul\`es and A.-K. {Cheik Ahamed}.
\newblock Alinea: An advanced linear algebra library for massively parallel
  computations on graphics processing units.
\newblock {\em Int. J. High Perform. Comput. Appl.}, 29(3):284--310, 2015.

\bibitem{Magoules2017b}
F.~Magoul\`es and G.~Gbikpi-Benissan.
\newblock {JACK}: An asynchronous communication kernel library for iterative
  algorithms.
\newblock {\em J. Supercomput.}, 73(8):3468--3487, 2017.

\bibitem{Magoules2018e}
F.~Magoul\`es and G.~Gbikpi-Benissan.
\newblock Distributed convergence detection based on global residual error
  under asynchronous iterations.
\newblock {\em IEEE Trans. Parallel Distrib. Syst.}, 29(4):819--829, 2018.

\bibitem{Magoules2018b}
F.~Magoul\`es and G.~Gbikpi-Benissan.
\newblock {JACK2}: An {MPI}-based communication library with non-blocking
  synchronization for asynchronous iterations.
\newblock {\em Adv. Eng. Softw.}, 119:116--133, 2018.

\bibitem{Magoules2018}
F.~Magoul\`es, G.~Gbikpi-Benissan, and Q.~Zou.
\newblock Asynchronous iterations of {P}arareal algorithm for option pricing
  models.
\newblock {\em Mathematics}, 6(4):1--18, 2018.

\bibitem{Magoules2015e}
F.~Magoul\`es, F.-X. Roux, and G.~Houzeaux.
\newblock {\em Parallel Scientific Computing}.
\newblock Wiley-ISTE, 2015.

\bibitem{Magoules2018c}
F.~Magoul\`es and C.~Venet.
\newblock Asynchronous iterative sub-structuring methods.
\newblock {\em Math. Comput. Simul.}, 145:34--49, 2018.

\bibitem{Yuan2006}
Y.-X. Yuan.
\newblock A new stepsize for the steepest descent method.
\newblock {\em J. Comput. Math.}, 24(2):149--156, 2006.

\bibitem{Zou2018}
Q.~Zou and F.~Magoul\`es.
\newblock A new cyclic gradient method adapted to large-scale linear systems.
\newblock In {\em Proceedings of the 17th International Symposium on
  Distributed Computing and Applications to Business, Engineering and Science},
  pages 196--199, Wuxi, China, 2018. IEEE.

\bibitem{ZouX1}
Q.~Zou and F.~Magoul\`es.
\newblock Fast gradient methods with alignment for symmetric linear systems
  without using {C}auchy step.
\newblock {\em preprint available at arXiv:1909.01479}, 2019.

\bibitem{ZouX2}
Q.~Zou and F.~Magoul\`es.
\newblock Parameter estimation in the {H}ermitian and skew-{H}ermitian
  splitting method using gradient iterations.
\newblock {\em preprint available at arXiv:1909.01481}, 2019.

\bibitem{ZouX4}
Q.~Zou and F.~Magoul\`es.
\newblock Recent developments in iterative methods for reducing
  synchronization.
\newblock {\em preprint available at arXiv:1912.00816}, 2019.

\bibitem{ZouX3}
Q.~Zou and F.~Magoul\`es.
\newblock Reducing the effect of global synchronization in delayed gradient
  methods for symmetric linear systems.
\newblock submitted.

\end{thebibliography}
\bibliographystyle{abbrv}

\end{document}